\documentclass[12 pt]{article}
\usepackage[american]{babel}
\usepackage{amsfonts}
\usepackage{latexsym}
\usepackage{amsmath}
\usepackage{amssymb}
\usepackage{amsthm}
\usepackage{fullpage}
\newtheorem{theorem}{Theorem}[section]
\newtheorem{lemma}[theorem]{Lemma}

\newtheorem{proposition}[theorem]{Proposition}
\newtheorem{corollary}[theorem]{Corollary}
\newtheorem{problem}[theorem]{Problem}

\newtheorem{conjecture}[theorem]{Conjecture}
\newtheorem{thmy}{Theorem}
\newenvironment{oldtheorem}{\stepcounter{thm}\begin{thmy}}{\end{thmy}}
\theoremstyle{definition}

\theoremstyle{remark}
\newtheorem*{note*}{Note}



\newcommand{\spana}{\textnormal{span}}
\begin{document}
\small

\title{\bf On the shape of a convex body with respect to its second projection body}

\medskip

\author {Christos Saroglou}

\date{August, 2014}

\maketitle

\begin{abstract}
\footnotesize We prove results relative to the problem of finding sharp bounds for the affine invariant $P(K)=V(\Pi K)/V^{d-1}(K)$. Namely,
we prove that if $K$ is a 3-dimensional zonoid of volume 1,
then its second projection body $\Pi^2K$ is contained in $8K$, while if $K$ is any symmetric
3-dimensional convex body of volume 1, then $\Pi^2K$ contains $6K$. Both inclusions are sharp. Consequences of these results include a stronger version
of a reverse isoperimetric inequality for 3-dimensional zonoids-established by the author in a previous work, a reduction for the 3-dimensional
Petty conjecture to another isoperimetric problem and the best known
lower bound up to date for $P(K)$ in 3 dimensions. As byproduct of our methods, we establish an almost optimal lower bound
for high-dimensional bodies of revolution.
\end{abstract}

\section{Introduction}
$\hspace{1em}$Let $K$ be a convex body in $\mathbb{R}^d$, that is convex, compact and with non-empty interior. Denote the volume of $K$ by $V_d(K)=V(K)$
and by $h_K$ its support function. The support function of $K$ is defined as
$$h_K(x)=\max_{y\in K}|\langle x,y\rangle|\ , \ x\in\mathbb{R}^d \ .$$ The support function of $K$ is convex and positively homogeneous. Moreover,
if $L$ is a different convex body with $L\supseteq K$, then $h_L\geq h_K$ and $h_K\neq h_L$, therefore $K$ is characterized by its support function. On
the other hand, any convex and positively homogeneous function is the support function of a (unique) convex body. The support function is, furthermore,
additive with respect to the Minkowski addition: $h_{K+M}=h_K+h_M$, where $M$ is a convex body and $K+M:=\{x+y\ |\ x\in K,\ y\in M\}$.

Define the \emph{projection body} $\Pi K$ of $K$
by its support function
$$h_{\Pi K}(u)=V_{d-1}\big(K  |  u^{\bot}\big)=\frac{1}{2}\int_{S^{d-1}}|\langle u,y\rangle | \ dS_K(y) \ , \ u\in S^{d-1} \ .$$
Here, $K  |  u^{\bot}$ denotes the orthogonal projection of $K$ on the subspace which is orthogonal to $u$ and $S(K,\cdot)$ stands for the surface area
measure of $K$, viewed as a measure on $S^{d-1}$. This is defined as:
$$S_K(\Omega)=V\Big( \bigl\{x\in \textnormal{bd}(K):\exists \ u\in \Omega,  \textnormal{ so that }
u\textnormal{ is a normal unit vector for } K \textnormal{ at } x\big\}\Big)\ . $$
It is clear that the projection body of a polytope is a \emph{zonotope}, i.e, the Minkowski sum of a finite number of line segments and in general the
projection body of a convex body is a \emph{zonoid}, i.e. a limit of zonotopes, in the sense of the Haussdorff metric. Conversely, it can be proven
that all zonoids are, up to translation, projection bodies of convex bodies. We refer to  \cite{Sch1} \cite{G} \cite{Ko} \cite{KRZ}
\cite{KRZ2} for more
information about support functions and
projection bodies.

C. M. Petty \cite{PE2} proved that $\Pi K$ is an affine equivariant, in particular the following expression$$P(K):=\frac{V(\Pi K)}{V(K)^{d-1}}$$
is invariant under invertible affine maps. Therefore, and since $P(K)$ is continuous with respect to the Haussdorff metric, classical arguments ensure
the existence of its extremal values; these are both finite and positive. The determination of these extremals is a difficult and challenging problem in
convex and affine geometry. It should be noted that the problem is interesting only for $d\geq 3$. In the plane (see e.g. \cite{Sa2}),
the minimizers are all centrally symmetric convex figures and the maximizers are precisely the triangles.
Before we state our main results, we would like to say a few words about the history of this problem.

Petty \cite{PE} conjectured that $$P(K)\geq \omega_{d-1}^d\omega_d^{2-d} \ ,$$
with equality if and only if $K$ is an ellipsoid, where $\omega_d$ is the volume of the $d$-dimensional unit ball $B_2^d$. The Petty conjecture, if true,
would be a powerful tool in the study of various functionals on convex bodies, since it would imply a whole family of isoperimetric inequalities (see
\cite{Lut1} \cite{Lut2} \cite{Lut3}, \cite{Lut4}), including the classical isoperimetric inequality and Petty projection inequality (see below).
Indicative of the difficulty
of the problem, is the fact that absolutely no
partial results are known to be true towards its solution. For instance, other major conjectures of variational nature such as the simplex conjecture
(for the Sylvester problem) or the Mahler conjecture, have been confirmed in special natural classes of convex bodies
(i.e. that contain naturally the conjectured optimizer) and/or are known to hold locally
(in a certain sense). We refer e.g. to \cite{B-F}, \cite{re1}, \cite{re2}, \cite{Stra}, \cite{NPRZ}, \cite{KR}, \cite{CCG}, \cite{Ra}.
No such results are valid for Petty's conjecture. We should mention here that Petty's conjecture holds for bodies with significantly larger
minimal surface area than the ellipsoid's, as follows from \cite{GP}, however these bodies are far from being ellipsoids, with respect to the
Banach-Mazur distance. Another
discouraging fact, proven by the author \cite{Sa2}, is that in general the \emph{Steiner symmetrization} fails for this problem.
Let us however mention the two most important
results that are known to be true:\\
(i) Petty's conjecture is equivalent (see \cite{Lut1}) to the following conjecture:
If $K$, $L$ are convex bodies of volume 1,
then the quantity
$$\int_{S^{n-1}}\int_{S^{n-1}}|\langle x, y\rangle |dS_K(x)dS_L(y)$$ is minimal when $K$ and $L$ are balls. \\
(ii) $P(K)\geq P(\Pi K)$, with equality if and only
if $K$ is homothetic to $\Pi^2K:=
\Pi(\Pi(K))$. This first appears in Schneider \cite{Sch3} and is called ``the class reduction technique'' (see section \ref{Section-Auxilliary-facts}).
This has been exploited in the study various other functionals (see e.g. \cite{Lut0} \cite{Lut2} \cite{Lut3} \cite{LYZ2}). It should be remarked here
that the previous result shows that if $K$ is a solution to Petty's problem, then $K$ is homothetic to $\Pi^2K$.
Therefore, any solution to the Petty problem has to be a zonoid (and in particular centrally symmetric). Ellipsoids have this property,
however not only these. Weil showed that the only polytopes
with this property are cartesian products of symmetric convex bodies of dimension 1 or 2. We will call such bodies "Weil bodies``. Note that the
three-dimensional Weil bodies are exactly the symmetric cylinders.
Our main results are very much
related to the class reduction technique; see Theorems \ref{theorem upper bound} and \ref{theorem lower bound 3-d} below.

The problem of the determination of the best upper bound for $P(K)$ is known as the Schneider problem. It was conjectured by Schneider \cite{Sch2} that
\begin{equation}
P(K)\leq 2^d\label{Schneider conjecture}
\end{equation}
in the centrally symmetric case, with equality precisely for the Weil-bodies. This conjecture was disproved by Brannen \cite{Br1}, who conjectured \cite{Br1}
\cite{Br2} that the centrally symmetric body of maximal volume contained in a simplex whose centroid is at the origin (see also \cite{Sa3})
is the solution in the symmetric case and the simplex
in the general case. The restriction of the problem to smaller classes of convex bodies also makes sense. For example, it was conjectured by Brannen that
that the
Weil-bodies are indeed the only maximizers for $P(K)$ in the class of zonoids. This was proven in three dimensions by Saroglou \cite{Sa2}
as of the upper bound;
however equality
occurs in other cases than the cylinders as well. To be more specific, the following was proven:

\begin{oldtheorem}\label{oldthm}
If $K$ is a 3-dimensional zonoid, then $P(K)\leq 2^3$, with equality if and only if $K$ is the sum
of five line segments or the sum of a cylinder and a line segment.
\end{oldtheorem}
The complicated cases of equality in Theorem \ref{oldthm}
is a rather discouraging fact in the attempt of extending this in higher dimensions.

We should mention here the following sharp inequalities, involving the polar of $\Pi K$,
$(\Pi K)^{\circ}=\{x\in \mathbb{R}^d|\langle x, y \rangle\leq 1, \ \forall y\in \Pi K\}$:
\begin{equation}
\frac{(2d)!}{d^d(d!)^2}\leq V\big((\Pi K)^{\circ}\big)V(K)^{d-1}\leq \big(\omega_{d}/\omega_{d-1}\big)^d\ ,
\label{petty-projection-ineq&zhang-ineq}
\end{equation}
with equality on the right if and only if $K$ is an ellipsoid and on the left if and only if $K$ is a simplex. The right inequality is the
Petty projection
inequality \cite{PE} (see also \cite{SM} for a different approach and \cite{LYZ2} for extension
to the $L^p$-projection bodies) and the left inequality is the Zhang inequality \cite{ZH} (see also \cite{G-Z}).
Combining (\ref{petty-projection-ineq&zhang-ineq})
with the Blaschke-Santal\'{o} \cite{Bl} \cite{Santalo} \cite{Me-Pa}
inequality and its reverse exact for zonoids  \cite{re1} \cite{re2},
$$\frac{4^d}{d!}\leq V(K)V(K^{\circ})\leq \omega_d^2\ ,$$
we obtain
\begin{equation}
C_2^{-d}P(B_2^d)\leq(C_1^{-1}/\sqrt{d})^d\leq\frac{4^d}{d!}\big(\omega_{d}/\omega_{d-1}\big)^{-d}\leq P(K)\leq \omega_d^2\frac{d^d(d!)^2}{(2d)!}\leq
C_1^d\leq C_2^dP(S_d) \ ,
\label{asymptotic-form}
\end{equation}
where $S_d$ is a $d$-dimensional simplex and $1<C_1<C_2$ are absolute constants. As far as we know, the estimate $P(K)\geq27/6$ that follows from
(\ref{asymptotic-form}) for 3-dimensional convex bodies was the best known lower bound up to date, in dimension 3. Other functionals involving
projection bodies were studied in \cite{LYZ} \cite{He-Le-Li}.

Our aim is to study the following
\begin{problem}\label{problem} If $\cal{C}$ is a certain class of convex bodies (we are particularly interested in the
class of symmetric
convex bodies and the class of zonoids) in $\mathbb{R}^d$, what are the smallest possible numbers $a,b>0$, such that
$$b^{-1}K\subseteq \Pi^2K\subseteq aK \ ,$$
for all $K\in{\cal{C}}$, of volume 1?
\end{problem}
It turns out that this problem is closely connected with the problem of finding the extremals for $P(K)$; this will be clear in Section 2.

Assume, now that $K$ contains the origin in its interior. It is convenient to work with the following functionals:
$$M(K):=\max_{x\in S^{d-1}}\frac{h_{\Pi^2K}(x)}{h_K(x)V(K)^{d-2}}$$
$$m(K):=\min_{x\in S^{d-1}}\frac{h_{\Pi^2K}(x)}{h_K(x)V(K)^{d-2}}\ .$$
The advantage of these quantities is that they both are centroaffine invariants (see next section). Let us, now, state our main results (Theorems
\ref{theorem upper bound} and \ref{theorem lower bound 3-d}):

\begin{theorem}\label{theorem upper bound}
Let $K$ be a zonoid in $\mathbb{R}^3$. Then,
\begin{equation}
M(K)\leq 2^3.\label{upper-bound-ineq}
\end{equation}
Equality holds if and only if $K$ is a centrally symmetric cylinder.
\end{theorem}
It is natural to conjecture the following:
\begin{conjecture}\label{conjecture upper bound}
Let $K$ be a zonoid in $\mathbb{R}^d$. Then, $\Pi^2 K\subseteq 2^dV(K)^{{d-2}}K$, with equality if and only if $K$ is a Weil-convex body.
\end{conjecture}
\begin{theorem}\label{theorem lower bound 3-d}
Let $K$ be a centrally symmetric convex body in $\mathbb{R}^3$. Then, for all $x\in S^2$,
$$\frac{h_{\Pi^2K}(x)}{h_K(x)V(K)}\geq m(K)\geq 6 \ ,$$
with equality if and only if $K$ is a symmetric double cone with circular base and $x$ is parallel to its main axis.
\end{theorem}
In other words, for zonoids, $\Pi^2K\subseteq 8V(K)K$ with equality only for cylinders and for general symmetric convex bodies,
$\Pi^2K\supseteq 6V(K)K$ and the boundary of $\Pi^2K$ can only touch the boundary of $6V(K)K$ in some direction $x$, if $K$ is a
symmetric double cone with circular base and $x$ is parallel to its main axis.

Actually, a stronger statement than Theorem \ref{theorem lower bound 3-d} will be shown to be true (Corollary \ref{lemma-last}),
which will lead to the following:

\begin{theorem}\label{3-d2revolution}
Let $K$ be a symmetric convex body in $\mathbb{R}^3$. The following is true
$$P(K)\geq \max_{x\in S^2}\Bigg\{4\left(\displaystyle\int_{-h_K(x)}^{h_K(x)} V_{2}\left( K\cap (sx+x^{\bot})\right)^{\frac{1}{2}}ds\right)^{2}
\Bigg/{ h_K(x)V(K)}\Bigg\}
\ .$$
\end{theorem}
We conjecture the following:
\begin{conjecture}\label{conjecture-3d-petty}
The quantity $$\max_{x\in S^2}\Bigg\{\left(\displaystyle\int_{-h_K(x)}^{h_K(x)} V_{2}\left( K\cap (sx+x^{\bot})\right)^{\frac{1}{2}}ds\right)^{2}
\Bigg/{ h_K(x)V(K)}\Bigg\}$$
is minimized precisely for (3-dimensional) ellipsoids.
\end{conjecture}
\begin{conjecture}\label{larger-dimensional-version-ofmain-thm-2}
Let $K$ be a $d$-dimensional convex body. Then,
$$m(K)\geq\max_{x\in S^{d-1}}\omega_{d-2}^{d-1}\omega_{d-1}^{3-d}\Bigg\{\left(\displaystyle\int_{-h_K(x)}^{h_K(x)} V_{d-1}\left( K\cap (sx+x^{\bot})\right)^{\frac{d-2}{d-1}}ds\right)^{d-1}
\Bigg/{ h_K(x)V(K)^{d-2}}\Bigg\}\ .$$
\end{conjecture}

Our method will also provide a much better estimate for bodies of revolution in $d$-dimensions
(compare with the general known bound $P(K)\geq c^d P(B_2^d)$, for some $c<1$), when $d$ is large enough.

\begin{theorem}\label{lower bound bodies of revolution}
Let $K$ be a convex body of revolution. Then,
$$P(K)> 2\frac{d^{d-2}}{(d-1)^{d-1}}\omega_{d-2}^{d-1}\omega_{d-1}^{3-d}= 2\frac{c_d}{\sqrt{d}}P(B_2^d)\ ,$$
where $c_d$ is a constant that depends only on $d$ and tends to $\sqrt{e/2\pi}$, as $d\rightarrow \infty$.
\end{theorem}
We believe that Problem \ref{problem} and Theorems \ref{theorem upper bound}, \ref{theorem lower bound 3-d} are of independend interest. We would like, however,
to summarize in brief the contribution of our results in the study of the extremals of $P(K)$.\\
Theorem \ref{theorem upper bound}:
\begin{itemize}
 \item Strenghtens significantly Theorem \ref{oldthm} (see Proposition \ref{2ndpb implies Brannen} below).
 \item Provides hope for extension in higher dimensions, since the equality cases in 3 dimensions are natural and expected, unlike Theorem \ref{oldthm},
 where the equality cases are much more complicated.
\end{itemize}
Our approach for the lower bound:
\begin{itemize}
 \item Improves the best known lower bound for $P(K)$ in 3 dimensions: $P(K)\geq 6$ instead of $P(K)\geq 27/6$ (combine Theorem
\ref{theorem lower bound 3-d} with Proposition \ref{conjecture-m(K)-implies-better-bound-for-P(K)} below).
 \item Reduces Petty's Conjecture in three dimensions (via Theorem \ref{3-d2revolution}) to Conjecture \ref{conjecture-3d-petty}.
 \item Provides hope for an almost optimal lower bound for $P(K)$ in high dimensions; this bound actually holds for bodies of revolution (Theorem
 \ref{lower bound bodies of revolution}).
\end{itemize}

\section{Auxiliary facts}\label{Section-Auxilliary-facts}
\hspace*{1em}Let us first check that $M(K)$ and $m(K)$ are invariant under linear maps.
The exponent $d-2$ in the denominators is chosen so that the ratios are invariant under
dilations, thus we only have to consider volume-preserving linear maps. Let us for instance prove that $M(K)$ is
invariant under $SL_n$-maps. For $T\in SL_n$, we have
\begin{eqnarray*}
M(TK)&=&\max_{x\in S^{d-1}}\frac{h_{\Pi^2(TK)}(x)}{h_{TK}(x)V(TK)^{d-2}}\\
&=&\max_{x\in S^{d-1}}\frac{h_{T\Pi^2K}(x)}{h_{TK}(x)V(K)^{d-2}}\\
&=&\max_{x\in S^{d-1}}\frac{h_{\Pi^2K}\big(T^tx\big)}{h_{K}\big(T^tx\big)V(K)^{d-2}}\\
&=&\max_{x\in S^{d-1}}\frac{h_{\Pi^2K}\big(T^tx/|T^tx|\big)}{h_{K}\big(T^tx/|T^tx|\big)V(K)^{d-2}}\\
&=&\max_{y\in S^{d-1}}\frac{h_{\Pi^2K}(y)}{h_{K}(y)V(K)^{d-2}}=M(K)\ .\\
\end{eqnarray*}

Next, we would like to describe the class reduction technique (promised in the introduction), as presented in \cite{Sch1}.
Let $K$, $L$ be convex bodies in $\mathbb{R}^d$. One can use Fubini's Theorem to obtain
\begin{equation}
 V(\Pi L, \Pi K)=V(\Pi^2K,L)\ .
\label{class-reduction-Fubini}
\end{equation}
Recall the definition of mixed volumes, their integral representation and the classical Minkowski inequality (see e.g. \cite{Sch1} or \cite{G}):
$$V(K,L):=V(K,L\dots,L):=\lim_{t\rightarrow 0^+}\frac{V(L+tK)-V(L)}{t}=\frac{1}{d}\int_{s^{n-1}}h_K(u)dS_L(u)\geq V(K)^{1/d}V(L)^{(d-1)/d}  .$$
Thus by (\ref{class-reduction-Fubini}) and  the Minkowski inequality, we get
\begin{eqnarray*}
V(\Pi K)&=&V(\Pi K, \Pi K)
= V(\Pi^2 K,K)
\geq V(\Pi^2K)^{1/d}V(K)^{(d-1)/d}\ ,
\end{eqnarray*}
with equality if and only if $K$ is homothetic to $\Pi^2 K$, which gives
$P(\Pi K)\geq P(K) ,$ with equality if and only if $K$ is homothetic to $\Pi^2 K$.
\begin{proposition}\label{2ndpb-P(K)-proposition}
Let $\mathcal{C}$ be a class of convex bodies in $\mathbb{R}^d$, which is topologically closed and closed under the action of the operator $\Pi$. Then,
\begin{equation}
\min_{K\in \mathcal{C}}P(K)=\min_{K\in \mathcal{C}}M(K) \ .\label{2ndpb-P(K)-equation}
\end{equation}
In particular, Petty's conjecture holds if and only if $M(K)$ is minimal when $K$ is an ellipsoid.
\end{proposition}
Proof. Take $K\in\mathcal{C}$, for which $M(K)$ is minimal. Then, by (\ref{class-reduction-Fubini}), we have:
\begin{eqnarray*}
V(\Pi K)&=&V(\Pi^2K, K)\\
&\leq& V\left(M(K)V(K)K,K\right)^{d-2}\\
&=&M(K)V(K)^{d-1}\ ,
\end{eqnarray*}
thus $\min_{K'\in \mathcal{C}}P(K')\leq P(K)\leq M(K)$. Assume, now, that $K\in\mathcal{C}$ is such that $P(K)$ is minimal.
Then, since $P(\Pi K)\leq P(K)$ and $\Pi K\in \mathcal{C}$, it follows that $P(\Pi K)=P(K)$,
which shows that $\Pi^2 K$ is homothetic to $K$. Thus, $\Pi^2K=\lambda K$, for some $\lambda>0$.
We may clearly assume that $V(K)=1$. It follows that $V(\Pi^2K)=\lambda^dV(K)=\lambda^d$.
Clearly, $M(K)=\lambda$.
Therefore,
\begin{eqnarray*}
\lambda^d&=&V(\Pi^2K)\\
&=&\frac{V(\Pi^2K)}{V(\Pi K)^{d-1}}\left(\frac{V(\Pi K)}{V(K)^{d-1}}\right)^{d-1}\\
&=&P(\Pi K)P(K)^{d-1}=P(K)^d \ ,
\end{eqnarray*}
thus $\min_{K'\in \mathcal{C}}M(K')\leq M(K)=\lambda=P(K)$. $\Box$
\begin{proposition}\label{conjecture-m(K)-implies-better-bound-for-P(K)}
\begin{flalign}
&(i)\ P(K)\geq m(K) \ ,\ \textnormal{for any convex body }K\in\mathbb{R}^n.&\nonumber \\
&(ii) \textnormal{ If conjecture } \ref{larger-dimensional-version-ofmain-thm-2}
\textnormal{ holds, then }P(K)\geq 2d^{d-2}(d-1)^{-(d-1)}\omega_{d-2}^{d-1}\omega_{d-1}^{3-d}>\big(c/\sqrt{d}\big)P\big(B_2^d\big)\ . &\nonumber
\label{P(K)>m(K)}
\end{flalign}

\end{proposition}
Proof. Assertion (i) follows again by (\ref{class-reduction-Fubini}), with $L=K$ and (ii) follows immediately from (i) and Lemma \ref{lemma-bodies-of-revolution} in
the end of this note. $\Box$\\
\\

Theorem \ref{3-d2revolution} follows immediately from Proposition \ref{conjecture-m(K)-implies-better-bound-for-P(K)} (i) and Corollary \ref{lemma-last}
(see below). Combining (\ref{class-reduction-Fubini}) and Theorem \ref{theorem upper bound}, it follows immediately:
\begin{proposition}\label{2ndpb implies Brannen}
If Conjecture \ref{conjecture upper bound} is true in $\mathbb{R}^d$, then for every zonoid $K$ and for every convex body $L$ in $\mathbb{R}^d$,
$$V(\Pi L,\Pi K)\leq 2^dV(K)^{d-2}V(K,L)\ .$$
In particular,
$$V(\Pi K)\leq 2^dV(K)^{d-1}\ . $$
\end{proposition}

\section{Formulas about projection bodies and their projections}
Let $Z=\sum_{i=1}^n[-x_i,x_i]$ be a zonotope in
$\mathbb{R}^d$. The support function of $Z$ is given by
\begin{equation}
h_Z(x)=\sum_{i=1}^n|\langle x,x_i\rangle| \ , \ x\in \mathbb{R}^d \ .
\label{formula 0.5}
\end{equation}
The volume of $Z$ is given by (see \cite{W2} for
proof and extensions):
\begin{equation}
V(Z)=2^d\sum_{\{i_1,\dots ,i_d\}\subseteq [n]}|
\det(x_{i_1},\dots ,x_{i_d})|=\frac{2^d}{d!} \sum_{i_1,\dots ,i_d\in
[n]}|\det(x_{i_1},\dots ,x_{i_d})  | \ ,
\label{formula1}
\end{equation}
where $[n]:=\{1,\dots,n\}$.
Thus, if $F_1,\dots,F_n$ are the facets of a polytope
$K$ in $\mathbb{R}^d$ with corresponding outer normal unit
vectors $x_1,\dots,x_n$, by the definition of $\Pi K$ and
(\ref{formula1}), we have:
\begin{equation}
V(\Pi K)=\sum_{\{i_1,\dots,i_d\}\subseteq [n]}V(F_{i_1}) \dots
V(F_{i_d})\cdot | \det(x_{i_1},\dots,x_{i_d})  |\label{formula1.5}
\end{equation}
We will also need formulas about the projection of $Z$ and its projection body. It is clear by (\ref{formula1}) that
\begin{eqnarray}\label{formula2}
V_{d-1}\big(Z|x^{\bot}\big)
&=&\frac{2^{d-1}}{(d-1)!}\sum_{i_1,\dots ,i_{d-1}\in [n]}|\textnormal{det}_{(d-1)
\times (d-1)}(x_{i_1}|x^{\bot},\dots,x_{i_{d-1}}|x^{\bot})|\nonumber\\
&=&\frac{2^{d-1}}{(d-1)!}\sum_{i_1,\dots ,i_{d-1}\in [n]}|\det(x_{i_1},\dots,x_{i_{d-1}},x)|\ .
\end{eqnarray}
As noted in \cite{Sa2}, one can show that when the $x_i$'s are in general position, then the
facets of $\Pi Z$ are up to translation exactly the parallelepipeds spanned by the vectors $x_{i_1},\dots x_{i_{d-1}}$,
where $i_1,\dots,
i_{d-1}\in [n]$. Projecting onto $x^{\bot}$ and using (\ref{formula1.5}), one easily obtains
\begin{eqnarray}\label{formula3}
(\Pi Z)|x^{\bot}
&=&2^{d-1}\sum_{1\leq i_1<\dots<i_{d-1}\leq n}[-(x_{i_1}\wedge\dots\wedge x_{i_{d-1}})|x^{\bot},
(x_{i_1}\wedge\dots\wedge x_{i_{d-1}})|x^{\bot}]\nonumber\\
&=&\frac{2^{d-1}}{(d-1)!}\sum_{i_1,\dots,i_{d-1}\in [n]}[-(x_{i_1}\wedge\dots\wedge x_{i_{d-1}})|x^{\bot},
(x_{i_1}\wedge\dots\wedge x_{i_{d-1}})|x^{\bot}] \ .
\end{eqnarray}
Combining (\ref{formula1}) and (\ref{formula3}), we immediately get:
\begin{eqnarray}\label{formula4}
V_{d-1}\big((\Pi Z)\ |\ x^{\bot}\big)=\left(\frac{2^{d-1}}{(d-1)!}\right)^d\sum_{i_1,\dots,i_{(d-1)^2}}
|\det(x_{i_1}\wedge\dots\wedge x_{i_{d-1}},\dots,x_{i_{(d-2)(d-1)+1}}\wedge\dots\wedge x_{i_{(d-1)^2}},x)|\ .
\end{eqnarray}
Finally, take $K$ to be a simplicial polytope in $\mathbb{R}^d$. Set $Vert(K)$ and $F(K)$ for the set of vertices and the set of facets of $K$
respectively. For $F\in F(K)$ and $v_1,\dots,v_{d}\in Vert(F)$, the outer normal vector on $F$ with respect to $K$, whose length equals the
$(d-1)$-dimensional volume of $F$ is exactly (up to a change of sign) $[1/(d-1)!](v_1-v_d)\wedge\dots \wedge (v_{d-1}-v_d)$. Therefore,
\begin{equation}
\Pi K=\frac{1}{2(d-1)!}\sum_{F\in F(K),\ v_1,\dots v_d\in Vert(F)}[-(v_1-v_d)\wedge\dots \wedge
(v_{d-1}-v_d),(v_1-v_d)\wedge\dots \wedge (v_{d-1}-v_d)] \ .
\label{formula5}
\end{equation}
\section{Upper inclusion}
The goal of this section is to establish Theorem \ref{theorem upper bound}. The proof will be a modification of the proof of Theorem \ref{oldthm}.
In the same spirit as in \cite{Sa2}, define the quantities:
$$S(x_1,\dots,x_{(d-1)^2},x)=\sum_{i_1,\dots,i_{(d-1)^2}\in [(d-1)^2]}|\det(x_{i_1},\dots x_{i_{d}})\dots
\det(x_{i_{(d-3)d+1}} ,\dots,x_{i_{d(d-2)}})|\cdot|\langle x_{i_{(d-1)^2}},x\rangle|$$
$$T(x_1,\dots,x_{(d-1)^2},x)=\sum_{i_1,\dots,i_{(d-1)^2}\in [(d-1)^2]}|\det(x_{i_1}\wedge\dots\wedge x_{i_{d-1}},\dots,
x_{i_{(d-2)(d-1)+1}}\wedge\dots\wedge x_{i_{(d-1)^2}})| \ ,$$
where $x_{i_1},\dots,x_{i_{(d-1)^2}},x\in \mathbb{R}^d$. It is clear that both functions $S$ and $T$ are
convex and positively homogeneous in each one of their arguments. Moreover, as it follows by the proof of the following lemma,
$S(x_1,\dots,x_{(d-1)^2},x)=0$ if and only if $T(x_1,\dots,x_{(d-1)^2},x)=0$.
\begin{lemma}\label{T-S-equivalence with inclusion-lemma}
The following are equivalent:\\
\begin{flalign}
&(i)\ \frac{h_{\Pi^2Z}(x)}{h_Z(x)V(Z)^{d-2}}\leq 2^d\ ,
\textnormal{ for every zonoid }Z\textnormal{ in }\mathbb{R}^d\textnormal{ and every }x\in S^{d-1} \ .& \nonumber \\
&(ii) \ T(x_1,\dots,x_{(d-1)^2},x)\leq \frac{[(d-1)!]^2}{d^{d-2}}S(x_1,\dots,x_{(d-1)^2},x),
\textnormal{ for all } x_{i_1},\dots,x_{i_{(d-1)^2}},x\in \mathbb{R}^d\ . &
\label{T-S-equivalence with inclusion-ineq}
\end{flalign}

\end{lemma}
Proof. Let $Z=\sum_{i=1}^n[-x_i,x_i]$ be a zonotope in $\mathbb{R}^d$ and $x\in S^{d-1}$. By (\ref{formula4}) one has
$$h_{\Pi^2Z}(x)=V\left((\Pi Z)|x^{\bot}\right)_{d-1}=\frac{2^{(d-1)d}}{[(d-1)!]^2}\sum_{i_1,\dots,i_{(d-1)^2}\in [n]}T(x_1,\dots,x_{(d-1)^2},x)\ .$$
Also, by (\ref{formula 0.5}), (\ref{formula1}), it can be easily proven that
$$h_Z(x)V(Z)^{d-2}=\frac{2^{d(d-2)}}{(d!)^{d-2}}\sum_{i_1,\dots,i_{(d-1)^2}\in [n]}T(x_1,\dots,x_{(d-1)^2},x)\ .$$
Therefore, by (\ref{formula4}),
\begin{equation}
\frac{h_{\Pi^2Z}(x)}{h_Z(x)V(Z)^{d-2}}
=\frac{2^dd^{d-2}}{[(d-1)!]^2}
\frac{\sum_{i_1,\dots,i_{(d-1)^2}\in [n]}T(x_1,\dots,x_{(d-1)^2},x)}{\sum_{i_1,\dots,i_{(d-1)^2}\in [n]}S(x_1,\dots,x_{(d-1)^2},x)}
\label{T-S-equivalence with inclusion-identity}
\end{equation}
Now, by (\ref{T-S-equivalence with inclusion-identity}), (ii) immediately implies (i). For the converse, take $n=(d-1)^2$. Leting
$|x_1|\rightarrow 0$, it follows by the homogeneity of $T$ and $S$ that terms of the form $T(x_1,x_1,\dots,x_{(d-1)^2},x)$ and
$S(x_1,x_1,\dots,x_{(d-1)^2},x)$ vanish in
(\ref{T-S-equivalence with inclusion-identity}) (and also all terms that correspond to permutations of the $x_i$'). Doing the same for the rest of the
$x_i$'s and using the symmetry of $T$ and $S$, assumption (i) and again (\ref{T-S-equivalence with inclusion-identity}), we can easily obtain
(\ref{T-S-equivalence with inclusion-ineq}).
$\Box$
\\

The rest of the section will be devoted in the proof of (\ref{T-S-equivalence with inclusion-ineq}) in dimension 3.
\begin{lemma}\label{convex/affine}
Let $f, \ g:[a,b]\rightarrow \mathbb{R}$ be functions such that $f,g$ are positive in $(a,b)$, $f$ is affine and $g$ is convex.
Assume, furthermore that the side limits $\lim_{t\rightarrow a^+}\frac{g(t)}{f(t)}$, $\lim_{t\rightarrow b^-}\frac{g(t)}{f(t)}$
are finite. Then, for every $t\in (a,b)$,
$$\frac{g(t)}{f(t)}\leq\max\Big\{\lim_{t\rightarrow a^+}\frac{g(t)}{f(t)}, \lim_{t\rightarrow b^-}\frac{g(t)}{f(t)}\Big\} \ .$$
Finally, if $f(a)=0$ (resp. $f(b)=0$), then for every $t\in (a,b)$,
$$\frac{g(t)}{f(t)}\leq \lim_{t\rightarrow b^-}\frac{g(t)}{f(t)}\
\textnormal{  (resp. }\frac{g(t)}{f(t)}\leq\lim_{t\rightarrow a^+}\frac{g(t)}{f(t)}) \ .$$
\end{lemma}
Proof. For $\lambda\in (0,1)$,
\begin{eqnarray}
\frac{g(\lambda a+(1-\lambda)b)}{f(\lambda a+(1-\lambda)b)}&=&
\frac{g(\lambda a+(1-\lambda)b)}{\lambda f(a)+(1-\lambda)f(b)}\nonumber\label{formula in g/f}\\
&\leq&\frac{\lambda g(a)+(1-\lambda)g(b)}{\lambda f(a)+(1-\lambda)f(b)}\\
&\leq&\max\Big\{\lim_{t\rightarrow a^+}\frac{g(t)}{f(t)}, \lim_{t\rightarrow b^-}\frac{g(t)}{f(t)}\Big\}\ .\nonumber
\end{eqnarray}
Moreover, if $f(a)=0$ (the case $f(b)=0$ is exactly the same; note also that by assumption the case $f(a)=f(b)=0$ cannot occur),
then $g(a)=0$, so by (\ref{formula in g/f}),
$$\frac{g(\lambda a+(1-\lambda)b)}{f(\lambda a+(1-\lambda)b)}\leq
\frac{0+(1-\lambda)g(b)}{0+(1-\lambda)f(b)}=\frac{g(b)}{f(b)}=\lim_{t\rightarrow b^-}\frac{g(t)}{f(t)}\ .\ \Box$$

We will also need some three-dimensional geometric lemmas.

\begin{lemma}\label{T=4/3S}
Let $x_1,x_2,x_3,x_4,x\in\mathbb{R}^3$. If two of the $x_i$'s are parallel, then
$$T(x_1,x_2,x_3,x_4,x)=\frac{4}{3}S(x_1,x_2,x_3,x_4,x) \ .$$
\end{lemma}
Proof. Since the function $(x_1,x_2,x_3,x_4)\mapsto(T/S)(x_1,x_2,x_3,x_4,x)$ is symmetric and 0-homogeneous, we may assume that $x_3=x_4$.
Using the identity $(x_1\wedge x_3)\wedge(x_2\wedge x_3)=\det(x_1,x_2,x_3)\cdot x_3$,
it follows that $$|\det(x_1\wedge x_3,x_2\wedge x_3, x)|=|\det(x_1,x_2,x_3)\langle x,x_3\rangle|\ .$$ One can check that
$T(x_1,x_2,x_3,x_3,x)=2^3|\det(x_1\wedge x_3,x_2\wedge x_3, x)|=8|\det(x_1,x_2,x_3)\langle x,x_3\rangle|$ and also $S(x_1,x_2,x_3,x_3,x)
=3!|\det(x_1,x_2,x_3)\langle x,x_3\rangle|. \ \Box$
\\

The proof of the next lemma can be easily understood by drawing a figure and we omit it.
\begin{lemma}\label{noproof}
Let $x_1,x_2,x_3,x_4,x\in \mathbb{R}^3$. Consider the planes $E_1:=\spana\{x_1,x_2\}$, $E_2:=\spana\{x_2,x_3\}$,
$E_3:=\spana\{x_1,x_3\}$, $E_4:=x^{\bot}$. Take $\nu\in\mathbb{R}^3$ and $t_1<0<t_2$ and suppose that one of the following holds:
\begin{flalign*}
&(i)\ x_4+t\nu\not\in E_j,\ j=1,2,3,4,\textnormal{ for all }t\in(t_1,t_2)\ . &\\
&(ii)\ x_4,\nu\in E_i,\textnormal{ for some } i\in[4],\textnormal{ and }x_4+t\nu\not\in E_j, \ i\in[4] \setminus\{i\} . &
\end{flalign*}
Then, the function $[t_1,t_2]\ni t\mapsto S(x_1,x_2,x_3,x_4+t\nu,x)\in \mathbb{R}$ is affine.
\end{lemma}
\begin{lemma}\label{main-techical-lemma}
Let $x_1,x_2,x_3,x_4,E_1,E_2,E_3,E_4$ be as in Lemma \ref{noproof}. Then, there exists a vector $x_4'$, such that:
\begin{flalign*}
&i)\textnormal{ There exist }i,j\in [4],i\neq  j,\textnormal{ such that }x_4'\in E_i\cap E_j\ .&\\
&ii)\ S(x_1,x_2,x_3,x_4,x)\neq 0 \textnormal{ and }\frac{T}{S}(x_1,x_2,x_3,x_4,x)\leq \frac{T}{S}(x_1,x_2,x_3,x_4',x) \ . &
\end{flalign*}

\end{lemma}
Proof. Set for simplicity $(T/S)(y_1,y_2,y_3,y_4,y):=0$, for any vectors $y_1,y_2,y_3,y_4,y\in \mathbb{R}^3$, with
$S(y_1,y_2,y_3,y_4,y)=0$. Suppose first that $x_4\in E_i$, for some $i\in [4]$. We may assume that
$x_4\not\in E_j$, for $j\in [4]\setminus\{i\}$. Otherwise we are done. It is obvious that the vector $x_4$
is contained in a convex angle spanned by the lines of the form $E_k\cap E_i$, $E_l\cap E_i$, $k,l\in [4]\setminus\{i\}$,
$k\neq l$ and also this angle is disjoint with $E_j\setminus \{0\}$, where $j=[4]\setminus \{i,k,l\}$. Take any vector $\nu\in E_i\setminus\{x_4,0\}$. Then, there exist $t_1<0<t_2$ (or $t_2<0<t_1$), such that
$x_4+t_1\nu\in E_k$ and $x_4+t_2\nu\in E_l$. According to Lemma \ref{noproof} (ii), the function
$t\mapsto S(x_1,x_2,x_3,x_4+t\nu,x)$ is affine in $[t_1,t_2]$, so by Lemma \ref{convex/affine},
$$(T/S)(x_1,x_2,x_3,x_4,x)=(T/S)(x_1,x_2,x_3,x_4+0\nu,x)\leq \max_{m=1,2}(T/S)(x_1,x_2,x_3,x_4+t_m\nu,x)\ .$$
Note that the last quantity cannot be zero, otherwise $S(x_1,x_2,x_3,x_4+t_m\nu,x)$ would be zero, $m=1,2$ which would force $S(x_1,x_2,x_3,x_4,x)$ to be zero.
Now, we may assume without loss of generality that
$(T/S)(x_1,x_2,x_3,x_4+t_1\nu,x)\geq (T/S)(x_1,x_2,x_3,x_4,x)$. Then, $x_4':=x_4+t_1\nu$ is as required.

It remains to deal with the case in which $x_4\not\in E_1,E_2,E_3,E_4$.
It is clear that there exist $t_1<0<t_2$ and $\nu\in\mathbb{R}^3$, such that $x_4+t_1\nu\in E_k$, $x_4+t_2\nu\in E_l$, for some $k,l\in[4]$, $k\neq l$ and
$x_4+t\nu\not\in E_1,E_2,E_3,E_4$, for all $t\in [t_1,t_2]$.
Therefore, by Lemma \ref{noproof} (i), the function
$t\mapsto S(x_1,x_2,x_3,x_4+t\nu,x)$ is again affine in $[t_1,t_2]$.
Using Lemma \ref{convex/affine} as before, we conclude that $(T/S)(x_1,x_2,x_3,x_4,x)\leq \max_{m=1,2}(T/S)(x_1,x_2,x_3,x_4+t_m\nu,x)$.
Since the last quantity cannot be zero, it is clear that we now fall in the previous special case where $x_4\in E_i$, for some $i\in[4]$ and our assertion is proved. $\Box$
\\
\\
Proof of Theorem \ref{theorem upper bound}:\\
By Lemma \ref{T-S-equivalence with inclusion-lemma}, it suffices to prove (\ref{T-S-equivalence with inclusion-ineq}) in 3 dimensions.
Let $x_1,x_2,x_3,x_4,x\in\mathbb{R}^3$, such that $S(x_1,x_2,x_3,x_4,x)\neq 0$. With the notation of Lemma
\ref{noproof}, we may assume by Lemma \ref{main-techical-lemma}
that $x_4\in E_1\cap E_j$, where $j\in [4]$. If $j\neq 4$, then $E_1\cap E_j=\spana\{x_k\}$, where $k\in\{1,2,3\}$.
Thus, $x_4$ is parallel to one of $x_1,x_2,x_3$, which by Lemma \ref{T=4/3S} gives $(T/S)(x_1,x_2,x_3,x_4,x)=4/3$.
Hence, we may assume that $x_4\in E_4=x^{\bot}$. Applying the same argument in $x_3$ instead of $x_4$, we may assume that $x_3\in x^{\bot}$ as well.
It remains to apply Lemma \ref{main-techical-lemma} one last time; $x_2$ can be replaced
by a vector $x_2'$, so that $x_2'$ is contained in a plane spanned by the $x_i$'s
and in a different plane spanned by the $x_i$'s or in $x^{\bot}$ and at the same time $S(x_1,x_2',x_3,x_4,x)\neq 0$.
However, if $x_2'\in x^{\bot}$, one can check that $S(x_1,x_2',x_3,x_4,x)= 0$, which shows that $x_2'$ is necessarily parallel
to $x_1$ or $x_3$ or $x_4$. This, by Lemma \ref{T=4/3S}, ends the proof of (\ref{upper-bound-ineq}).

We are left with the equality cases. First take $x_1,x_2,x_3,x_4\in \mathbb{R}^3$, not all coplanar, and assume that they have the propertry that
$T(x_1,x_2,x_3,x_4,x)=(4/3)S(x_1,x_2,x_3,x_4,x)$, for all $x\in S^2$. Consider the zonotope $W=\sum_{i=1}^4[-x_i,x_i]$.
Then, (\ref{T-S-equivalence with inclusion-identity}) together
with Lemma \ref{T=4/3S} yields
$$\frac{h_{\Pi^2W}(x)}{h_W(x)V(W)}=2^3\frac{3}{4}\frac{T(x_1,x_2,x_3,x_4,x)+\sum_{\{i_1,i_2,i_3,i_4\}\subsetneq [4]}T(x_{i_1},x_{i_2},x{i_3},x{i_4},x)}
{S(x_1,x_2,x_3,x_4,x)+\sum_{\{i_1,i_2,i_3,i_4\}\subsetneq [4]}S(x_{i_1},x_{i_2},x{i_3},x{i_4},x)}=2^3\ , $$
for all $x\in S^2$. Therefore, $W$ is a polytope which is homothetic to its second projection body, so by Weil's result mentioned in the introduction,
$W$ is a cylinder.

Now, let $Z$ be any zonoid in $\mathbb{R}^3$. Its support function is given by $h_Z(x)=\int_{S^2}|\langle x,y\rangle|d\mu (y)$,
for some even measure $\mu$ in $S^2$. Then, (\ref{T-S-equivalence with inclusion-identity}) immediately implies
$$\frac{h_{\Pi^2Z}(x)}{h_Z(x)V(Z)}=2^3\frac{3}{4}\frac{\int_{S^2}\int_{S^2}\int_{S^2}\int_{S^2}T(x_1,x_2,x_3,x_4,x)d\mu(x_1)d\mu(x_2)d\mu(x_3)d\mu(x_4)}
{\int_{S^2}\int_{S^2}\int_{S^2}\int_{S^2}S(x_1,x_2,x_3,x_4,x)d\mu(x_1)d\mu(x_2)d\mu(x_3)d\mu(x_4)}\  , \ x\in S^2 \ .$$
Note that $S$ and $T$ are continuous functions. Thus, if $h_{\Pi^2Z}(x)/(V(Z)h_Z(x))=8$,
for all $x\in S^2$, then $(T/S)(x_1,x_2,x_3,x_4,x)=4/3$, for all $x\in S^2$
and for all $x_1,x_2,x_3,x_4$ from the support of $\mu$. This shows that for any $x_1,x_2,x_3,x_4$ in the support of $\mu$,
the zonotope $W=\sum_{i=1}^4[-x_i,x_i]$
is a cylinder, which leads to the fact that $Z$ is a cylinder itself. $\Box$
\section{Symmetrization and lower bounds}
Let $K$ be convex body in $\mathbb{R}^d$ and $\nu$ be a unit vector. The Steiner symmetrization $S_{\nu}K$ of $K$
along the direction $\nu$ is defined as the convex body, whose intersection with every line $l$ parallel to $\nu$ is
symmetric with respect to $\nu^{\bot}$ and has the same length as $l\cap K$. Closely related is the so called Schwartz
symmetrization $T_{\nu}K$ of $K$, which is defined as the convex body whose intersection with with every hyperplane $H$
orthogonal to $\nu$ is a ball whose center belongs to the line $\mathbb{R}\nu$ and has the same $(d-1)$-dimensional volume
as $K\cap H$. Both $S_{\nu}K$ and $T_{\nu}K$ are convex bodies, as follows by the Brunn-Minkowski Theorem.
Furthermore, it is true there exists a sequence of directions, so that
if Steiner symmetrization is applied along these directions successively, the resulting sequence of convex bodies will converge to a ball.
Similarly, there exists a sequence of
directions from $\nu^{\bot}$, so that the corresponding sequence of convex
bodies will converge to $T_{\nu}K$. It follows by the Fubini theorem that $S_{\nu}(K)$ and $T_{\nu}(K)$ preserve the volume of $K$.

Let us define the Steiner symmetrization of $K$ in a more convenient way. One can write:
$$K=\{x+t\nu:x\in K|{\nu^{\bot}},\ f(x)\leq t\leq g(x)\} \ ,$$ where $g,-f:K|{\nu^{\bot}}\rightarrow \mathbb{R}^d$ are concave functions.
Then, $$S_{\nu}K=\{x+t\nu:x\in K|{\nu^{\bot}},\ -w(x)\leq t\leq w(x)\}\ ,$$ where
$w=(g-f)/2$. Set $u:=(f+g)/2$. Define also the function $u_K:K|\nu^{\bot}\rightarrow \mathbb{R}$, with
$u_K(x):=u(x|\nu^{\bot})$, $x\in K|\nu^{\bot}$. It is clear that
\begin{equation*}
K=\{x+u_K(x)\nu:x\in S_{\nu}K\}\ ,\ K^{\nu}=\{x-u_K(x)\nu:x\in S_{\nu}K\}\ ,\label{K-K^{nu}-def}
\end{equation*}
where $K^{\nu}$ is the reflection of $K$ with respect to the hyperplane $\nu^{\bot}$.

We will need the following lemma, which can follow from \cite[Lemma 4]{Sa1}. Actually, \cite[Lemma 4]{Sa1} refers to a more general class of transformations of convex bodies,
the so-called shadow systems; these were systematically studied by Campi and Gronchi/Campi, Colesanti and Gronchi ; see e.g. \cite{CCG} \cite{CG}.
\begin{lemma}\label{STeiner-symmetrization-saroglou} Assume that $K$ is a simplicial polytope in $\mathbb{R}^d$.
Then, $S_{\nu}(K)$ is also a simplicial polytope. Suppose, in addition, that for every vertex $v$ of $K$, there exists another vertex $v'$ of $K$, such that the
line segment $[v,v']$ is parallel to $\nu$. Then, the following is true: The facets of $K$ are exactly the $(d-1)$-dimensional simplices
$\textnormal{conv} \{v_i+u_K(v_i)\nu:v_i\in V(F),\ i=1\dots d\}$, where $F\in F(S_{\nu}(K))$.
\end{lemma}
Notice that the class of (symmetric) polytopes satisfying the assumptions of Lemma \ref{STeiner-symmetrization-saroglou}
is dense in the class of (symmetric) convex bodies in $\mathbb{R}^d$.
Thus, one may use the previous lemma at will, when trying to prove inequalities involving Steiner symmetrization.
\begin{lemma}\label{lemma-lower-main-technical}
Let $\nu\in S^{d-1}$, $H$ be a 2-dimensional subspace of $\mathbb{R}^d$, containing $\nu$ and $K$ be a convex body in $\mathbb{R}^d$.
Then, $V_2\left([\Pi(S_{\nu}K)]\ | \ H\right)\leq V_2([\Pi K]\ | \ H)$.
\end{lemma}
Proof. We may assume that the assumptions of Lemma \ref{STeiner-symmetrization-saroglou} are satisfied (in particular $S_{\nu}K$ is a simplicial polytope).
Acoording to (\ref{formula5}),
$$\Pi(S_{\nu}K)=a_d\sum_{{F\in F(S_{\nu}K)}\atop{v_1,\dots,v_d\in Vert(F)}}[\pm(v_1-v_d)\wedge\dots\wedge (v_{d-1}-v_d)] \ ,$$
where $[\pm w]:=[-w,w]$ and $a_d$ is a constant which depends only on the dimension $d$. It is clear by Lemma  \ref{STeiner-symmetrization-saroglou}
and (\ref{formula5}), that
$$\Pi K=a_d\sum_{{F\in F(S_{\nu}K)}\atop{v_1,\dots,v_d\in Vert(F)}}[\pm ((v_1-u_K(v_1)\nu)-(v_d-u_K(v_d)\nu))\wedge\dots\wedge
((v_{d-1}-u_K(v_{d-1})\nu-(v_d-u_K(v_d)\nu))] \ ,$$
$$\Pi( K^{\nu})=a_d\sum_{{F\in F(S_{\nu}K)}\atop{v_1,\dots,v_d\in Vert(F)}}[\pm ((v_1+u_K(v_1)\nu)-(v_d+u_K(v_d)\nu))\wedge\dots\wedge
((v_{d-1}+u_K(v_{d-1})\nu-(v_d+u_K(v_d)\nu))] \ .$$
Let $F$ be a facet of $S_{\nu}K$ and $v_1,\dots,v_d$ be the vertices of $F$. Set
$$U_F:=(v_1-v_d)\wedge\dots\wedge (v_{d-1}-v_d) \ ,$$
$$N_F:=\sum_{i=1}^d(v_1-v_d)\wedge\dots\wedge(v_{i-1}-v_d)\wedge(u_K(v_d)-u_K(v_i))\nu\wedge (v_{i+1}-v_d)\wedge\dots\wedge (v_{d-1}-v_d) \ .$$
It is, then, clear that
$$\Pi(S_{\nu}K)=a_d\sum_{{F\in F(S_{\nu}K)}}[-U_F,U_F] \ ,$$
$$\Pi K=a_d\sum_{{F\in F(S_{\nu}K)}}[-U_F-N_F,U_F+N_F] \ ,
\ \Pi (K^{\nu})=a_d\sum_{{F\in F(S_{\nu}K)}}[-U_F+N_F,U_F-N_F] \ .$$
Now, without loss of generality, we may take $H=\spana \{e_1,e_2\}$, $\nu=e_2$. Notice that for $F\in F(S_{\nu}K)$, $N_F$ is always orthogonal to $\nu=e_2$. Thus,
$N_F|H=\langle N_F,e_1\rangle e_1$. Set $U_F^i=\langle U_F,e_i \rangle$ and $N_F^i=\langle N_F,e_i \rangle$, $i=1,2$.
We have shown that
$$(\Pi K)\ |\ H=a_d\sum_{{F\in F(S_{\nu}K)}}[-U_F^1e_1-U_F^2e_2-N_F^1e_1,U_F^1e_1+U_F^2e_2+N_F^1e_1] \ ,$$
therefore using (\ref{formula1}) we obtain
\begin{equation}
V_2\big( (\Pi K)\ |\ H\big)=
b_d\sum_{{F_1,F_2\in F(S_{\nu}K)}}\left|\det\left( (U_{F_1}^1+N_{F_1}^1,U_{F_1}^2),(U_{F_2}^1+N_{F_2}^1,U_{F_2}^2)\right) \right| \ ,
\label{111}
\end{equation}
where $b_d$ is another constant that depends on $d$. Similarly, we get
\begin{equation}
V_2\big( [\Pi (K^{\nu})]\ |\ H\big)=
b_d\sum_{{F_1,F_2\in F(S_{\nu}K)}}\left|\det\left( (U_{F_1}^1-N_{F_1}^1,U_{F_1}^2),(U_{F_2}^1-N_{F_2}^1,U_{F_2}^2)\right) \right| \ ,
\label{222}
\end{equation}
\begin{equation}
V_2\big( [\Pi(S_{\nu}K)]\ | \ H\big)=b_d\sum_{{F_1,F_2\in F(S_{\nu}K)}}\left|\det\left( (U_{F_1}^1,U_{F_1}^2),(U_{F_2}^1,U_{F_2}^2)\right) \right| \ .
\label{333}
\end{equation}
Observe that $V_2((\Pi K)\ |\ H)=V_2([\Pi (K^{\nu})]\ |\ H )$, thus by (\ref{111}), (\ref{222}) and (\ref{333}) we clearly have
$$V_2\big((\Pi K)\ |\ H\big)=\frac{V_2\big(\Pi K\ |\ H\big)+V_2\big([\Pi (K^{\nu})]\ |\ H\big)}{2}\geq V_2\big([\Pi(S_{\nu}K)]\ |\ H\big)\ ,$$
ending the proof. $\Box$
\\

Assume, now, that $K$ is a convex body in $\mathbb{R}^3$. For $x,\nu\in S^2$, $\nu\in x^{\bot}$ and since,
$V_2((\Pi K)\ |\ x^{\bot})=h_{\Pi^2K}(x)$, it follows that
$h_{\Pi^2(S_{\nu}K)}(x)\leq h_{\Pi^2K}(x)$. Since this is true for all $\nu\in x^{\bot}$,
it follows that $h_{\Pi^2(T_{x}K)}(x)\leq h_{\Pi^2K}(x)$. Note also that
by definition, $h_K(x)=h_{T_xK}(x)$. Hence, by Lemma \ref{lemma-lower-main-technical}, we have shown the following:
\begin{lemma}\label{lemma-schwartz-symmetrization}
Let $K$ be a convex body in $\mathbb{R}^3$ and $x\in S^2$. Then,
$$\frac{h_{\Pi^2(T_xK)}(x)}{h_{T_xK}(x)V(T_xK)}\leq \frac{h_{\Pi^2K}(x)}{h_K(x)V(K)} \ .$$
\end{lemma}

\begin{lemma}\label{lemma-bodies-of-revolution}
Let $K$ be a centrally symmetric convex body of revolution in $\mathbb{R}^d$ and $\mathbb{R}\nu$ be its axis of revolution 
for some $\nu\in S^{d-1}$. Then,
$$h_{\Pi^2K}(\nu)=\omega_{d-2}^{d-1}\omega_{d-1}^{3-d}\left(\int_{-h_K(\nu)}^{h_K(\nu)}V_{d-1}\left( K\cap (s\nu+\nu^{\bot})\right)^{\frac{d-2}{d-1}}ds
\right)^{d-1}\ .$$
\end{lemma}
Proof. Let $\mu$ be a direction orthogonal to $\nu$. One has the representation $K\cap \spana\{\nu,\mu\}=\{s\nu+t\mu:-h_K(\nu)\leq s\leq h_K(\nu),
\ -f(s)\leq t\leq f(s)\}$, where $f:[-h_K(\nu),h_K(\nu)]\rightarrow \mathbb{R}$
is a non-negative concave function. Clearly, for every $(d-1)$-dimensional subspace
$H$ of $\mathbb{R}^d$ that contains $\nu$, we have $V_{d-2}(K\cap H\cap (\nu^{\bot}+s\nu))=\omega_{d-2}f^{d-2}(s)$, $s\in [-h_K(\nu),h_K(\nu)]$, thus
\begin{equation}\label{identity 800}
V_{d-1}(K|H)=V_{d-1}(K\cap H)=\int_{-h_K(\nu)}^{h_K(\nu)}\omega_{d-2}f(s)^{d-2}ds \ .
\end{equation}

On the other hand, $V(K\cap (\nu^{\bot}+s\nu))_{d-1}=\omega_{d-1}f^{d-1}(s)$, so by (\ref{identity 800}) we have
\begin{equation}\label{identity 801}
V(K|H)_{d-1}=\omega_{d-2}\omega_{d-1}^{-\frac{d-2}{d-1}}\int_{-h_K(\nu)}^{h_K(\nu)}V_{d-1}\left( K\cap (\nu^{\bot}+s\nu)\right)^{\frac{d-2}{d-1}}ds\ .
\end{equation}
Now, since the projection of $\Pi K$ on the hyperplane $\nu^{\bot}$ is a ball of radius $h_{\Pi K}(x)$, $x\in \nu^{\bot}$, it follows that
$V_{d-1}((\Pi K)|\nu^{\bot})=\omega_{d-1}V_{d-1}(K|H)^{d-1}$, which by (\ref{identity 801}) proves our assertion. $\Box$\\

Combining the fact that the Schwartz symmetrization $T_{\nu}K$ of $K$
is a body of revolution, together with Lemmas \ref{lemma-schwartz-symmetrization} and
\ref{lemma-bodies-of-revolution}, we obtain the following.
\begin{corollary}\label{lemma-last}
Let $K$ be a symmetric convex body in $\mathbb{R}^3$ and $\nu\in S^2$. The following is true
$$\frac{h_{\Pi^2K}(\nu)}{h_K(\nu)V(K)}\geq 4\max_{x\in S^2}
\Bigg\{\left(\displaystyle\int_{-h_K(x)}^{h_K(x)} V_{2}\left( K\cap (sx+x^{\bot})\right)^{\frac{1}{2}}ds\right)^{2}
\Bigg/{ h_K(x)V(K)}\Bigg\}\ .$$
\end{corollary}
\begin{lemma}\label{use-of-berwald}
Let $K$ be a symmetric convex body of revolution in $\mathbb{R}^d$, whose axis of revolution is parallel to some unit vector $\nu$. Then,
$$\frac{h_{\Pi^2K}(\nu)}{h_K(\nu)V(K)^{d-2}}\geq  2\frac{d^{d-2}}{(d-1)^{d-1}}\omega_{d-2}^{d-1}\omega_{d-1}^{3-d}\ ,$$
with equality if and only if $K$ is a double cone whose main axis is parallel to $\nu$.
\end{lemma}
Proof.
By Lemma \ref{lemma-bodies-of-revolution} and the Fubini Theorem, we have:
\begin{eqnarray*}\label{berwald}
\frac{h_{\Pi^2K}(\nu)}{h_K(\nu)V(K)^{d-2}}&\geq&
\dfrac{\omega_{d-2}^{d-1}\omega_{d-1}^{3-d}\left(\displaystyle\int_{-h_K(\nu)}^{h_K(\nu)}V_{d-1}\left( K\cap (s\nu+\nu^{\bot})\right)^{\frac{d-2}{d-1}}
ds\right)^{d-1}}
{h_{K}(\nu)\left(\displaystyle\int_{-h_K(\nu)}^{h_K(\nu)}V_{d-1}\left( K\cap (s\nu+\nu^{\bot})\right)ds \right)^{d-2}} \\
&=&\omega_{d-2}^{d-1}\omega_{d-1}^{3-d}\dfrac{\left(\displaystyle\int_{-a}^{a}f^{d-1}(s)
ds\right)^{d-1}}{a\left( \displaystyle\int_{-a}^af^{d-1}(s)ds\right)^{d-2}}\ ,
\end{eqnarray*}
where $a=h_K(\nu)$ and $f(s)=V_{d-1}\left( K\cap (s\nu+\nu^{\bot})\right)^{1/d-1}$. Note that by the Brunn principle, the function $f:[-a,a]\rightarrow
\mathbb{R}$ is even and concave. The proof now follows immediately by the Berwald inequality \cite{Ber} (with $p=d-2$ and $q=d-1$):
If $f:[-a,a]\rightarrow \mathbb{R}_{+}$, is an even concave function and $0<p<q$, then
$$
\left(\frac{1}{2a}(1+p)\int_{-a}^af^p(s)ds\right)^{1/p}\geq  \left(\frac{1}{2a}(1+q)\int_{-a}^af^q(s)ds\right)^{1/q}\ ,
$$
with equality if and only if $f(s)=|f(0)-s(f(0)/a)|$. $\Box$
\\
\\
Proof for Theorem \ref{theorem lower bound 3-d}:\\
Immediate by Proposition \ref{conjecture-m(K)-implies-better-bound-for-P(K)} and Corollary \ref{lemma-last}. $\Box$
\\
\\
Proof of Theorem \ref{lower bound bodies of revolution}:\\
Since the class ${\cal R}$ of bodies of revolution, whose main axis is parallel to $\nu$, is closed under the action of the operator $\Pi$,
Proposition \ref{2ndpb-P(K)-proposition}
shows that $$\min_{K\in {\cal R}}P(K)=\min_{K\in {\cal R}}M(K)\geq \min_{K\in {\cal R}}\frac{h_{\Pi^2K}(\nu)}{V(K)^{d-2}h_K(\nu)} $$
and the proof follows by Lemma \ref{use-of-berwald}. $\Box$

\bigskip

\noindent \textsc{Ch.\ Saroglou}: Department of Mathematics,
Texas A$\&$M University, 77840 College Station, TX, USA.

\smallskip

\noindent \textit{E-mail:} \texttt{saroglou@math.tamu.edu \ \&\ christos.saroglou@gmail.com}

\end{document}